\documentclass[psfceqn,reqno]{amsart}
\usepackage{amsfonts}
\usepackage{amsthm}
\usepackage{amsmath}
\usepackage{amssymb}
\usepackage{color}
\usepackage{mathrsfs}
\usepackage{tikz}
\usetikzlibrary{arrows}
\usetikzlibrary{graphs}

\newtheorem{theorem}{Theorem}[section]
\newtheorem{corollary}[theorem]{Corollary}
\newtheorem{proposition}[theorem]{Proposition}
\newtheorem{lemma}[theorem]{Lemma}

\theoremstyle{definition}
\newtheorem{definition}[theorem]{Definition}
\newtheorem{remark}[theorem]{Remark}
\newtheorem{example}[theorem]{Example}

\def\R{{\mathbb R}}

\numberwithin{equation}{section}

    \thanks{\hspace{-1.66em} 2010 \emph{Mathematics Subject Classification.}
	Primary 91B50; Secondary 52HA99, 54C60, 54H25.}
    \thanks{\noindent $^\dagger$The first author was supported by a National Natural Science Foundation of 
   	China grant, No. NSFC-11271178. He would also like to thank the Department of Mathematical Sciences
   	at Auckland University of Technology for her hospitality during his visit in February 2019.\\
    $^\ddagger$ The corresponding author.}

\begin{document}
    \title[Existence of equilibrium $\cdots$]{Existence of equilibrium in an abstract\\ economy
    without $SS$-convexity}
    \author[J. Hou]{Ji-Cheng Hou$^\dagger$}
    \address{School of Science, Beijing Information Science and Technology University, Beijing 100192, 
    P. R. China}
    \email{e-mail:houjc163@163.com}
    \author[J. Cao]{Jiling Cao$^\ddagger$}
    \address{Department of Mathematical Sciences, School of Engineering, Computer and Mathematical Sciences, 
    Auckland University of Technology, Private Bag 92006, Auckland 1142, New Zealand}
    \email{jiling.cao@aut.ac.nz}
    \maketitle

    \begin{abstract}
    In this paper, we introduce the concept of $\Delta$-convexity in an abstract economy with all 
    choice sets being general topological spaces. We give a new generalization of the classical 
    Yannelis-Prabhakar equilibrium existence theorem in abstract economies by relaxing $SS$-convexity 
    to $\Delta$-convexity.
    \end{abstract}

    \section{Introduction} \label{sec:intro} 
    
    Throughout this paper, we use $2^S$ to denote the family of all subsets of a set $S$, and suppose all 
    topological spaces are Hausdorff. Let $A$ be a subset of a topological space $X$. We use $\overline{A}$ 
    to denote the closure of $A$ in $X$.  Let $X$ and $Y$ be two topological spaces. Let $T : X \rightarrow 
    2^Y$ be a correspondence. Then (i) for each $M \subseteq Y $, we denote by $T^{-1}(M)$ the set 
    $\{x\in X : T(x)\cap M \ne \emptyset\}$, and write $T^{-1}(y)$ for $T^{-1}(M)$ in the case that 
    $M =\{y \}$; (ii) $T$ is said to \emph{have open lower sections} if $T^{-1}(y)$ is open in $X$ for each 
    $y\in Y$; (iii) $\mbox{cl}T: X\mapsto 2^Y $ denotes the correspondence which assigns $\overline{T(x)}$ 
    to each point $x\in X$. If $A$ is a subset of a vector space, we denote by $\mbox{co}A$ the convex 
    hull of $A$. We use ${\R}$ to denote the set of all real numbers, ${\R}^{n+1}$ to denote the 
    $n+1$ dimensional Euclidean space, $\Delta_n$ to denote the standard $n$-dimensional simplex in 
    $\R^{n+1}$, and $e_i$ for $i=0, 1, \cdots, n$ to denote the standard base in $\Delta_n$. Let $I$ be 
    the set of agents that is either finite or infinite (even uncountable). Each agent $i$'s choice set 
    is a general topological space. Denote by $X$ the Cartesian product of all $X_i$'s equipped with the 
    product topology. For each $i\in I$,we denote by $x_{-i}$ the choices of all other agents rather 
    than agent $i$. Also denote by $X_{-I'} =\prod_{j\in I\setminus I'}X_j$ the Cartesian product of the 
    sets of choices of  agents $j\notin I'$, and we sample write $X_{-i}$ for the set $\{i\}$ consists 
    of a single point $i$.
    
    \medskip
    Following Borglin and Keiding \cite{Borglin:1976}, an abstract economy $\Gamma =(X_i, A_i, 
    P_i)_{i\in I}$ is defined as a family of ordered triples $(X_i, A_i, P_i)_{i\in I}$, where $I$ is 
    the set of all agents, $X_i$ is the choice set for agent $i$, $A_i:\prod_{i\in I}X_i\rightarrow 
    2^{X_i}$ is the constraint correspondence of agent $i$ , and  $P_i:\prod_{i\in I}X_i\rightarrow 
    2^{X_i}$ is the preference correspondence of agent $i$. An equilibrium for $\Gamma$ is an $x^*\in 
    X=\prod_{i\in I}X_i$ satisfying for each $i\in I$, $x^*_i\in \overline{A_i(x^*)}$ and $A_i(x^*)\cap 
    P_i(x^*)=\emptyset$. In \cite{Yannelis:1983} and \cite{Yannelis:1984}, Yannelis and Prabhakar  
    proved that, for countably many agents and countably many locally convex topological vector spaces, 
    an equilibrium for the abstract economy $\Gamma $ exists if for each $i\in I$,
    \begin{itemize}
    \item[(i)] $X_i$ is a nonempty, compact, convex, metrizable subset of a locally convex topological 
    vector space; 
    \item[(ii)] $A_i(x)$ is convex and nonempty for all $x\in X$; 
    \item[(iii)] the correspondence $\mbox{cl}A_i$ is upper semicontinuous;
    \item[(iv)] $A_i$ has open lower sections; 
    \item[(v)] $P_i$ has open lower sections; and
    \item[(vi)] $x_i\notin \mbox{co}P_i(x)$ for all $x\in X$.
    \end{itemize}
    Since then, this classical theorem has been extended in many ways, e.g., Tian \cite{Tian:1990}, 
    He and Yannelis \cite{He:1992}, Chowdhury et al. \cite{Chowdhury:2001}, Ding and Xia 
    \cite{Ding:2004}, Hou \cite{Hou:2008}, Kim and Yuan \cite{Kim:2001}, Shen \cite{Shen:2001}, 
    Tarafdar \cite{Tarafdar:1991}, Ionescu Tulcea  \cite{Tulcea:1988} and Yuan \cite{Yuan:2000}. 
    However, the results in all of these papers assume that the condition ``$x_i\notin \mbox{co}P_i(x)$ 
    for all $x\in X$" holds. This condition is called $SS$-convexity\footnote{Here, we use SS to stand 
    for Shafer and Sonnenschein, refer to \cite{Shafer:1975}.} in \cite{Tian:1994}.
    
    \medskip
    The purpose of this paper is to  give a new generalization of the classical Yannelis-Prabhakar 
    equilibrium existence theorem by relaxing the $SS$-convexity to a new notion called $\Delta$-convexity. 
    It should be emphasized that the method that we use is totally different from those methods given 
    in all results mentioned above.

    \section{The main result} \label{sec:notion} 
    
    In this section, we present our main result. First, we give the following simple lemma which 
    provides seemingly a new  approach on the investigation for the existence of equilibria in 
    abstract economies.
     
    \begin{lemma} \label{lem:new_approach}
    Let ${\Gamma}=(X_i,A_i,P_i)_{i\in I}$ be an abstract economy with $A_i(x)\ne \emptyset$ for 
    each $i\in I$ and each $x\in X$. Then ${\Gamma}$ has an equilibrium if and only if 
    \[
    \bigcap_{i\in I}\bigcap_{x_i\in X_i}\left[\left(F_i\cap \left(X\setminus P_i^{-1}(x_i)
    \right)\right) \cup \left(X\setminus A_i^{-1}(x_i)\right)\right]\ne \emptyset,
    \] 
    where $F_i=\left\{x\in X|x_i\in \overline{A_i(x)}\right\}$.
    \end{lemma}

    \begin{proof}
    \emph{Sufficiency.} Since 
    \[
    \bigcap_{i\in I}\bigcap_{x_i\in X_i}\left[\left(F_i\cap \left(X\setminus P_i^{-1}(x_i)
    \right)\right) \cup \left(X\setminus A_i^{-1}(x_i)\right)\right]\ne \emptyset,
    \] 
    we can pick up an element 
    \[
    x^*\in \bigcap_{i\in I}\bigcap_{x_i\in X_i}\left[\left(F_i\cap \left(X\setminus P_i^{-1}(x_i)
    \right)\right) \cup \left(X\setminus A_i^{-1}(x_i)\right)\right].
    \]
    However,
    \begin{eqnarray*}
    & & \bigcap_{i\in I}\bigcap_{x_i\in X_i}\left[\left(F_i\cap \left(X\setminus P_i^{-1}(x_i)
    \right)\right) \cup \left(X\setminus A_i^{-1}(x_i)\right)\right]\\
    &= & \bigcap_{i\in I}\bigcap_{x_i\in X_i} \left[\left(F_i\cup \left(X\setminus A_i^{-1}(x_i)
    \right)\right)\cap \left(\left(X\setminus P_i^{-1}(x_i)\right)\cup \left(X\setminus A_i^{-1}(x_i)
    \right)\right)\right].
    \end{eqnarray*}
    Therefore, 
    \[
    x^*\in \bigcap_{i\in I} \bigcap_{x_i\in X_i} \left(\left(X\setminus P_i^{-1}(x_i)\right) 
    \cup \left(X\setminus A_i^{-1}(x_i)\right)\right),
    \] 
    and thus $A_i(x^*)\bigcap P_i(x^*)=\emptyset $, for each $i\in I$. We now show that $x^*\in F_i$ 
    for each $i\in I$, which means that $x_{i^*}\in \overline{A_i(x^*)}$ for all $i\in I$. If not, 
    then there exists an $i^*\in I$ such that $x^*\notin F_{i*}$.  So, 
    \[
    x^*\in \bigcap_{x_{i^*}\in X_{i^*}} \left(X\setminus A_{i^*}^{-1}(x_{i^*})\right).
    \] 
    It follows that $A_{i^*}(x^*)=\emptyset $, which contradicts with the assumption that $A_i(x)\ne 
    \emptyset $ for all $i\in I$ and $x\in X$. Hence, $x^*$ is an equilibrium of $\Gamma $.
    
    \medskip
    {\it Necessity.}  If $x^*$ is an equilibrium of $\Gamma$, then 
    \begin{equation} \label{eqn1}
    x^*\in F_i \mbox{ and } A_i(x^*)\cap P_i(x^*)=\emptyset 
    \end{equation}
    for each $i\in I$.  We demonstrate that 
    \[
    x^*\in \bigcap_{i\in I}\bigcap_{x_i\in X_i}\left[\left(F_i\cap \left(X\setminus P_i^{-1}(x_i)\right)
    \right)\cup \left(X\setminus A_i^{-1}(x_i)\right)\right].
    \] 
    If not, then there exist an $i^*\in I$ and an $x_{i^*}\in X_i$ such that  
    \[
    x^*\notin \bigcap_{i\in I}\bigcap_{x_i\in X_i}\left[\left(F_i\cap \left(X\setminus P_i^{-1}(x_i)
    \right)\right)\cup \left(X\setminus A_i^{-1}(x_i)\right)\right].
    \] 
    Since $x^*\in F_{i^*}$, then we have $x^*\in P_{i^*}^{-1}(x_{i^*}))\cap A_{i^*}^{-1}(x_{i^*})$ 
    which further implies that $x_{i^*}\in P_{i^*}(x^*)\cap A_{i^*}(x^*)$, contradicting with
    \eqref{eqn1}. 
    \end{proof}
    
    In order to extend the Yannelis-Prabhakar equilibrium existence theorem to general topological 
    spaces, we need to make some preparation. First, we introduce the concept of $\Delta$-convexity
    in an abstract economy.
    
    \begin{definition} \label{def:delta_convex}
    An abstract economy $\Gamma =(X_i, A_i, P_i)_{i\in I}$ is called \emph{$\Delta$-convex}, if for 
    each $i\in I$  and each finite subset $\{x_{i0}, x_{i1}, \cdots, x_{in_i}\}$ of $X_i$, there exists 
    a continuous mapping $\varphi_{n_i}: \Delta_{n_i}\rightarrow X_i$ such that for any subset $J_i$ 
    of $\{0, 1, \cdots, n_i\}$,
    \begin{itemize}
    \item[(i)] if $x\in \bigcap_{j\in J_i} A_i^{-1}(x_{ij})$, then $\varphi_{n_i}(\Lambda_i)\in A_i(x)$ 
    for each 
    \[
    \Lambda_i\in \{(\lambda_{i0},\lambda_{i1}, \cdots,  \Lambda_{in_i})\in 
    \Delta_{n_i}: \lambda_{ij}=0, \mbox{ for all }j\notin J_i\};
    \] 
    and
    \item[(ii)] if $x\in \bigcap_{j\in J_i}P_i^{-1}(x_{ij})$, then $\varphi_{n_i}(\Lambda_i)\in P_i(x)$ 
    for each 
    \[
    \Lambda_i\in \{(\lambda_{i0},\lambda_{i1},\cdots, \lambda_{in_i})\in \Delta_{n_i}: \lambda_{ij}=0, 
    \mbox{ for all }j\notin J_i\}.
    \]
    \end{itemize}
    \end{definition}
    
    In order to demonstrate the rationales of Definition \ref{def:delta_convex}, we give the following 
    proposition.
    
    \begin{proposition} \label{prop:delta_convex}
    Let $\Gamma =(X_i, A_i, P_i)_{i\in I}$ be an abstract economy which satisfies that for each $i\in I$,
    \begin{itemize}
    \item[(i)] $X_i$ is a convex subset of a topological vector space; and
    
    \item[(ii)] $A_i(x)$ and $P_i(x)$ are convex for all $x\in X$.
    \end{itemize}
    Then $\Gamma =(X_i, A_i, P_i)_{i\in I}$ is $\Delta$-convex.
    \end{proposition}
    
    \begin{proof}
    For each $i\in I$ and each finite subset $\{x_{i0}, x_{i1}, \cdots, x_{in_i}\}\subseteq X_i$, 
    we define a mapping $\varphi_{n_i}:\Delta_{n_i}\rightarrow X_i$ by $\varphi_{n_i}(\Lambda_i)=\sum_{j=0}^{n_i}\lambda_{ij}x_{ij}$ for each $(\lambda_{i0},\cdots, 
    \lambda_{in_i})\in \Delta_{n_i}$. Obviously, $\varphi_{n_i}$ is continuous.  Let $J_i$ be 
    an arbitrary subset of $\{0, 1, \cdots, n_i\}$. Now, we show that (i) and (ii) in Definition 
    \ref{def:delta_convex} are satisfied. If $x\in \bigcap_{j\in J_i} A_i^{-1}(x_{ij})$, then $x_{ij}
    \in A_i^{-1}$ for all $j\in J_i$. Let 
    \[
    \Lambda_i\in \{(\lambda_{i0},\lambda_{i1}, \cdots, 
    \Lambda_{in_i})\in \Delta_{n_i}: \lambda_{ij}=0, \mbox{ for all }j\notin J_i\}.
    \] 
    Then, by the convexity of $A_i(x)$, one has 
    \[
    \varphi_{n_i}(\Lambda_i)=\sum_{j\in J_i}\lambda_{ij}x_{ij}\in A_i(x).
    \] 
    Using the same argument, we can verify (ii) in Definition \ref{def:delta_convex}. 
    \end{proof}
    
    \begin{remark}
    As shown in Example \ref{exam} later, the $\Delta$-convexity does not imply (ii) in Proposition 
    \ref{prop:delta_convex}.
    \end{remark}
    
    \begin{lemma}[\cite{Peleg:1967}] \label{lem:Peleg}
    Let $I_0=\{1,2,\cdots,n\}$ and $n_i$ be a natural number for each $i\in I_0$. For each $i\in I_0$, 
    let 
    \[
    F_i:\left\{e^i_j:j=0,1,2,\cdots, n_i\right\}\rightarrow \prod_{i\in I_0}\Delta_{n_i}
    \] 
    be a closed-valued set-valued map, where $\{e^i_j:j=0,1,2,\cdots, n_i\}$ denotes the standard 
    basis of $\Delta_{n_i}$. If for any $i\in I_0$ and any finite subset $S_i$ of $\{0,1,2,\cdots, 
    n_i\}$, one has 
    \[
    \prod_{i\in I_0}\mbox{co}\{e_j^i: j\in S_i\}\subseteq \bigcap_{i\in I_0}\bigcup_{j\in S_i}
    F_i\left(e^i_j\right).
    \] 
    Then $\bigcap_{i\in I_0}\bigcap_{j=0}^{n_i} F_i\left(e^i_j\right)\ne \emptyset$.
    \end{lemma}

    Now we state our main result.
    
    \begin{theorem} \label{thm:main}
    Let $\Gamma =(X_i, A_i, P_i)_{i\in I}$ be an abstract economy which satisfies that for each 
    $i\in I$,
    \begin{itemize}
    \item[(i)] $X_i$ is compact;
    
    \item[(ii)] $A_i(x)\ne \emptyset $ for all $x\in X$;
    
    \item[(iii)] the correspondence $\mbox{cl}A_i: X\rightarrow 2^{X_i}$ defined by $\mbox{cl}A_i(x)
    =\overline{A_i(x)}$ for all $x\in X$ is upper semi-continuous;
    
    \item[(iv)] $A_i$ has open lower section; 
    
    \item[(v)] $P_i$ has open lower section; and
    
    \item[(vi)] $x_i\notin P_i(x)$ for all $x\in X$, where $x_i$ is the projection of $x$ into $X_i$.
    \end{itemize}
    If $\Gamma$ is $\Delta$-convex, then it has an equilibrium.
    \end{theorem}

    \begin{proof}
    We use $K_i(x_i)$ to denote the set 
    \[
    \left[F_i\cap(X\setminus P_i^{-1}(x_i)\right]\cup[X\setminus A_i^{-1}(x_i)],
    \] 
    where 
    \[
    F_i=\left\{x=(x_i)_{i\in I}\in X: x_i\in \mbox{cl}A_i(x)\right\}.
    \]
    By Lemma \ref{lem:new_approach}, we only need to show that $\bigcap_{i\in I}\bigcap_{x_i\in X_i}
    K_i(x_i)\ne \emptyset$.
    
    \medskip
    \noindent 
    {\bf Claim 1.}  
    $\bigcap_{i\in I}\bigcap_{x_i\in X_i}\overline{K_i(x_i)}=\bigcap_{i\in I}\bigcap_{x_i\in X_i}K_i(x_i)$.
    
    \begin{proof}[Proof of Claim 1.]
    We only need to verify
    \[
    \bigcap_{i\in I}\bigcap_{x_i\in X_i}\overline{K_i(x_i)}\subseteq \bigcap_{i\in I}
    \bigcap_{x_i\in X_i}K_i(x_i).
    \] 
    Suppose the contrary. Then, there exist a 
    \begin{equation} \label{eqn:main}
    y\in \bigcap_{i\in I}\bigcap_{x_i\in X_i}\overline{K_i(x_i)},
    \end{equation} 
    an $i\in I$ and an $x_i\in X_i$ such that $y\notin K_i(x_i)$. So, $y\notin F_i\cap \left(X\setminus 
    P_i^{-1}(x_i)\right)$ and $y\notin X\setminus A_i^{-1}(x_i)$.
    
    \medskip
    \emph{Case 1.} $y\notin F_i$.
    By the definition of $F_i$, we have $(y,y_i)\notin Gr({\rm cl}A_i)$, where $Gr({\rm cl}A_i)$ is 
    the graph of ${\rm cl}A_i$ and $y_i$ is the projection of $y$ into $X_i$. By (iii) and the regularity 
    of $X\times X_i$, it is well known that $Gr({\rm cl}A_i)$ is closed in $X\times X_i$. Therefore, 
    there exists an open neighborhood $V$ of $(y,y_i)$ such that
    \begin{equation} \label{eqn:main2}
    V\cap Gr({\rm cl}A_i)= \emptyset.
    \end{equation}
    Without loss of generality, we may assume that $V=V(y)\times V(y_i)$ and 
    \[
    V(y)=\prod_{j\in J}V(y_j)\times \prod_{j\in J}X_j,
    \] 
    where $J$ is a finite subset of $I$ , $i\in J$, and $V(y_j)$ is an open neighborhood of $y_j$ for 
    all $j\in J$. Let $U(y)=V(y_i)\times\prod_{j\ne i}X_j$ and $W(y)=U(y)\cap V(y)$. We show 
    that $W(y)\cap F_i=\emptyset$. If not, we can pick up a point $z\in W(y)\cap F_i$. Then, $(z,z_i)
    \in Gr({\rm cl}A_i)$ and 
    \[
    z\in W(y)= U(y)\cap V(y)\subseteq \left(V(y_i)\times\prod_{j\ne i}X_j\right)\cap 
    \left(\prod_{j\in J}V(y_j)\times\prod_{j\notin J}X_j\right).
    \] 
    It follows that $(z,z_i)\in V$, which contradicts with \eqref{eqn:main2}. Thus, by \eqref{eqn:main} 
    and (iv), there exists an open neighborhood $Q(y)$ of $y$ such that $Q(y)\subseteq A^{-1}(x_i)$. 
    Hence, we have
    \[
    (Q(y)\cap W(y))\cap K_i(x_i)=\emptyset,
    \] 
    which again contradicts with \eqref{eqn:main2}.
    
    \medskip
    \emph{Case 2.} $y\in F_i$. Since $y\notin F_i\cap \left(X\setminus P_i^{-1}(x_i)\right)$ and 
    $y\notin X\setminus A_i^{-1}(x_i)$, we have $y\in P_i^{-1}(x_i)\cap A_i^{-1}(x_i)$. By (iv) 
    and (v), there exists an open neighborhood $Q(y)$ of $y$ such that $Q(y)\subseteq A^{-1}(x_i)
    \cap P_i^{-1}(x_i)$. Hence $Q(y)\cap K_i(x_i)=\emptyset$, which contradicts with \eqref{eqn:main}.
    This completes the proof of Claim 1.
    \end{proof}
    
    In order to complete the proof of this theorem, we need to show that 
    \[
    \bigcap_{i\in I}\bigcap_{x_i\in X_i}\overline{K_i(x_i)}\ne \emptyset.
    \] 
    Since $X$ is compact, we only need to show that the family 
    \[
    \left\{\overline{K_i(x_i)}: i\in I, x_i\in X_i \right\}
    \] 
    has the finite intersection property. Toward this end, let $I_0$ be an arbitrary finite subset 
    of $I$, and $B_i=\{b_{i1}, b_{i2}, \cdots, b_{in_i}\}$ be a finite subset of $X_i$ for each 
    $i\in I_0$. By the $\Delta$-convexity of $\Gamma$, for each $i\in I_0$, there exists a continuous 
    mapping $\varphi_{n_i}:\Delta_{n_i}\rightarrow X_i$ such that for any subset $J_i$ of $\{0, 1, 
    \cdots, n_i\}$,
    \begin{itemize}
    \item[(c1)] if $x\in \bigcap_{j\in J_i} A_i^{-1}(b_{ij})$, then $\varphi_{n_i}(\Lambda_i)\in 
    A_i(x)$ for any $\Lambda_i=(\lambda_{i0}, \cdots, \lambda_{in_i})\in \Delta_{n_i}$ with 
    $\lambda_{ij}=0\mbox{ for all }j\notin J_i$,
    
    \item[(c2)] if  $x\in \bigcap_{j\in J_i}(P_i^{-1}(b_{ij})$ , then $\varphi_{n_i}(\Lambda_i)\in 
    P_i(x)$ for any $\Lambda_i=(\lambda_{i0}, \cdots, \lambda_{in_i})\in \Delta_{n_i}$ with 
    $\lambda_{ij}=0 \mbox{ for all }j\notin J_i$.
    \end{itemize}
    Take an arbitrary point $x_{-I_0}^0\in X_{-I_0}$. Define a mapping $\Phi: \prod_{i\in I_0}
    \Delta_{n_i}\rightarrow X$ by
    \[
    \Phi((\Lambda_i)_{i\in I_0})=((\varphi_{n_i}(\Lambda_i))_{i\in I_0},x_{-I_0}^0)
    \] 
    for each $(\Lambda_i)_{i\in I_0}\in \prod_{i\in I_0}\Delta_{n_i}$. Obviously, $\Phi$ is a 
    continuous mapping from $\prod_{i\in I_0}\Delta_{n_i}$ into $X$. For each $i\in I_0$, we take 
    an arbitrary finite subset $S_i$ of $\{0, 1, \cdots, n_i\}$.
    
    \medskip
    \noindent 
    {\bf Claim 2.} $\prod_{i\in I_0}\mbox{co}\{e_j^i:j\in S_i\}\subseteq \bigcap_{i\in I_0}\bigcup_{j\in S_i}\Phi^{-1}\left(\overline{K_i(b_{ij})} \right).$
    
    \begin{proof}[Proof of Claim 2.] Suppose that $(\Lambda_i)_{i\in I_0}\notin \bigcap_{i\in I_0}
    \bigcup_{j\in S_i}\Phi^{-1}\left(\overline{K_i(b_{ij})}\right)$. Then, there exists an $i^*\in 
    I_0$ such that $(\Lambda_i)_{i\in I_0}\notin \Phi^{-1}\left(\overline{K_{i^*}(b_{i^*j})}\right)$ 
    for each $j\in S_{i^*}$. It follows that $\Phi((\Lambda_i)_{i\in I_0})\notin\overline{K_{i^*}
    	(b_{i^*j})}$ for each $j\in S_{i^*}$. Particularly, $\Phi((\Lambda_i)_{i\in I_0})\notin 
    K_{i^*}(b_{i^*j})$ for each $j\in S_{i^*}$. Let $x^*=\Phi((\Lambda_i)_{i\in I_0})$. Then $x^*_i=\varphi_{n_i}(\Lambda_i)$ for each $i\in I_0$, and $x^*_{-I_0}=x_{-I_0}^0$, by the 
    definition of $\Phi$.
    
    \medskip
    Since $x^*\notin K_{i^*}(b_{i^*j})$ for each $j\in S_{i^*}$, we have that
    \begin{equation} \label{eqn:main3}
    x^*\in \bigcap_{j\in S_{i^*}}A_{i^*}^{-1}(b_{i^*j})
    \end{equation}
    and
    \begin{equation} \label{eqn:main4}
    x^*\notin F_{i^*}\bigcap (X\setminus P_{i^*}^{-1}(b_{i^*j}))
    \end{equation} 
    for each $j\in S_{i^*}$. We show that $\Lambda_{i^*}\notin \mbox{co}\{e_j^{i^*}:j\in S_{i^*}\}$. 
    If not, then $\Lambda_{i^*}=\sum_{j\in S_{i^*}}\lambda_je_j^{i^*}$, where $\lambda_j\ge 0$ for 
    each $j\in S_{i^*}$ and $\sum_{j\in S_{i^*}}\lambda_j=1$. By (c1) and \eqref{eqn:main3}, $x^*_{i^*}=\varphi_{n_{i^*}}(\Lambda_{i^*})\in A_{i^*}(x^*)$, and thus $x^*\in F_{i^*}$. It 
    follows from \eqref{eqn:main4} that 
    \begin{equation} \label{eqn:main5}
    x^*\in \bigcap_{j\in S_{i^*}}P_{i^*}^{-1}(b_{i^*j})
    \end{equation} 
    Combing (c2) and \eqref{eqn:main5} together, we have
    \[
    x^*_{i^*}=\varphi_{n_{i^*}}(\Lambda_{i^*})\in P_{i^*}(x^*), 
    \] 
    which contradicts with the assumption that $x_i\notin P_i(x)$ for each $i\in I$ and $x\in X$. 
    This completes the proof of Claim 2.
    \end{proof}
    
    By using Lemma \ref{lem:Peleg} with $F_i(e^i_j)=\Phi^{-1}\left(\overline{K_i(b_{ij})}\right)$ 
    for each $i\in I_0$ and $j=0,1,\cdots, n_i$, we have $\bigcap_{i\in I_0} \bigcap_{j=0}^{n_i} \Phi^{-1}\left(\overline{K_i(b_{ij})}\right)\ne \emptyset$. We have completed the proof of
    the theorem.
    \end{proof}

    From Theorem \ref{thm:main} and Proposition \ref{prop:delta_convex}, we obtain immediately the 
    following corollary which improves the Yannelis-Prabhakar equilibrium existence theorem 
    mentioned previously, in the following three respects:
    \begin{enumerate}
    \item the local convexity of topological spaces for each $i\in I$ can be dropped;
    \item the condition that $X_i$ is metrizable for each $i\in I$ can be dropped; and
    \item the index set $I$ may be uncountable.
    \end{enumerate}
    
    \begin{corollary} \label{coro:main}
    Let $\Gamma =(X_i, A_i, P_i)_{i\in I}$ be an abstract economy which satisfies that for each 
    $i\in I$:
    \begin{itemize}
    \item[(i)] $X_i$ is a compact convex subset of a topological vector space;
    
    \item[(ii)] $A_i(x)$ is convex and nonempty for all $x\in X$;
    
    \item[(iii)] the correspondence ${\rm cl} A_i: X\rightarrow 2^{X_i}$, defined by ${\rm cl}A_i(x)
    =\overline{A_i(x)}$ for all $x\in X$, is upper semicontinuous;
    
    \item[(iv)] $A_i$ has open lower section; 
    
    \item[(v)] $P_i$ has open lower section; and
    
    \item[(vi)] $x_i\notin \mbox{co}P_i(x)$ for all $x\in X$, where $x_i$ is the projection of $x$ 
    into $X_i$.
    \end{itemize}
    Then $\Gamma$ has an equilibrium.
    \end{corollary}
    
    \begin{proof} 
    Instead of $P_i$ by ${\rm co}P_i$ defined as $({\rm co}P_i)(x)={\rm co}P_i(x)$ for all 
    $x\in X$. By Lemma 5.1 in \cite{Shafer:1975}, ${\rm co}P_i$ has open lower sections. According 
    to Proposition \ref{prop:delta_convex}, the abstract economy $\Gamma =(X_i, A_i, {\rm co}
    P_i)_{i\in I}$ satisfies all assumptions in Theorem \ref{thm:main}. It follows immediately that 
    $\Gamma =(X_i, A_i, P_i)_{i\in I}$ has an equilibrium. 
    \end{proof}

    \section{An example} \label{sec:quasi-cont} 
    
    In this section, we give an example of an abstract economy, in which all assumptions in Theorem 
    \ref{thm:main} hold. However, in this example, as the $SS$-convexity is not satisfied, not all 
    assumptions in Corollary \ref{coro:main} or the Yannelis-Prabhakar equilibrium existence theorem 
    hold.
    
    \begin{example} \label{exam}
    Consider an abstract economy $\Gamma=(X_i,A_i,P_i)_{i\in \{1,2\}}$ consisting of two agents, 
    where $X_1$ and $X_2$ are the closed intervals $[-1,1]$ and $[0,1]$, respectively, and the 
    restraint correspondences $A_i$ and preference correspondences $P_i$ are defined respectively 
    as follows:
    \begin{eqnarray*}
    	A_1(x_1,x_2)&=&\{1\},\\
    	A_2(x_1,x_2)&=&[0,1],\\
    	P_1(x_1,x_2)&=& \left(-1, -|x_1| \right)\cup \left(|x_1|, 1\right),\\
    	P_2(x_1,x_2)&=&\left\{
    	\begin{array}{ll}
    		\emptyset, &\ \mbox{if }|x_1|\ge x_2,\\
    		\left(0, |x_1|\right), &\ \mbox{if } |x_1|<x_2,
    	\end{array}
    	\right.
    \end{eqnarray*} 
    for each $(x_1,x_2)\in X=X_1\times X_2$. Obviously, $x_1\in \mbox{co}P_1(x)$ for each 
    \[
    x=(x_1, x_2)\in (-1,1)\times[0,1].
    \] 
    Therefore, the Yannelis and Prabhakar equilibrium existence theorem and Corollary \ref{coro:main}
    cannot be used.
    
    \medskip
    First, we verify the $\Delta$-convexity of $\Gamma$. For each $i\in \{1,2\}$ and each finite subset 
    $\{x_{i0}, x_{i1}, \cdots, x_{in_i}\}$ of $X_i$, we define mappings $\varphi_{n_i}:\Delta_{n_i}
    \rightarrow X_i$ by putting $\varphi_{n_i}(\Lambda_i)=\sum_{j=0}^{n_i}\lambda_{ij}|x_{ij}|$ for each $\Lambda_1=(\lambda_{10},\lambda_{11}, \cdots, \lambda_{1n_1})\in \Delta_{n_1}$. It is obvious that 
    $\varphi_{n_i}$ are continuous for all $i\in \{1,2\}$. Let $J_i$ be a finite subset of $\{i0, i1, 
    \cdots, in_i\}$. By setting $A_1(x)=A_2(x)=\{1\}$ , we do nothing for (i) of Definition 
    \ref{def:delta_convex}. If $x\in \bigcap_{j\in J_1}P_1^{-1}(x_{1j})$, then $x_{1j}\in P_1(x)$ for 
    all $j\in J_1$. Therefore, $|x_1|<|x_{1j}|<1$ for all $j\in J_i$. It follows that
    \[
    |x_1|<\varphi_{n_1}(\Lambda_1)=\sum_{j=0}^{n_1}\Lambda_{1n_1}|x_{ij}|<1,
    \] 
    which means $\varphi_{n_1}(\Lambda_1)\in P_1(x)$ for all 
    \[
    \Lambda_{1}\in \{(\lambda_{10}, \lambda_{11}, \cdots,\lambda_{1n_1})\in \Delta_{n_1}: \lambda_{1j} =0 
    \mbox{ for all }j\notin J_i\}.
    \] 
    If $x\in \bigcap_{j\in J_1}P_2^{-1}(x_{2j})$, then by the convexity of $P_2(x)$, we have
    $\varphi_{n_2}(\Lambda_2)\in P_1(x)$ for all 
    $\Lambda_{2}\in \{(\lambda_{20}, \lambda_{21}, \cdots,\lambda_{2n_2})\in \Delta_{n_2}:\lambda_{2j} =0 
    \mbox{ for all }j\notin J_2\}$.
    
    \medskip
    Now, we demonstrate that the assumptions (i)-(vi) in Theorem \ref{thm:main} are satisfied. We do 
    nothing for checking assumptions (i)-(iii) and (vi). Obviously, $A_i$ has open lower sections for 
    each $i=1,2$. We show that $P_i$ has also open lower sections for each $i=1,2$. Take an arbitrary 
    element $y_1\in X_1$ and an arbitrary element $x=(x_1,x_2)\in P_1^{-1}(y_1)$. Then, 
    \[
    y_1\in P_1(x)=(-1,-|x_1|)\cup (|x_1|,1),
    \] 
    i.e., $y_1>|x_1|$. Let $V=(-y_1,y_1)\times X_2$. It is easy to show that $V\subseteq P_1^{-1}(y_1)$.
    Take an arbitrary element $y_2\in X_2$ and an arbitrary element $x=(x_1,x_2)\in P_2^{-1}(y_2)$. Then 
    $|x_1|<x_2$ and $y_2\in (-|x_1|,|x_1|)$. Note $x_1\ne 0$. Without loss of generality, we support 
    $x_1>0$. Take an element $r\in (x_1,x_2)$  and put $V_1=(|y_2|,r)$ and $V_2=(r, 1]$. Obviously, 
    $x\in V_1\times V_2$. Take an arbitrary $x'=(x_1',x_2')\in V_1\times V_2 $. Then, 
    \[
    |y_2|<x_1'<r<x_2',
    \]
    which means $y_2\in P_2(x')=(-x_1',x_1')$, and thus $V_1\times V_2\subseteq P_2^{-1}(y_2)$.
    By Theorem \ref{thm:main}, $\Gamma$ has an equilibrium. \hfill $\Box$
    \end{example}
    
    From Theorem \ref{thm:main}, we can also obtain a new equilibrium existence theorem for qualitative 
    game \cite{Gale:1978}, by relaxing the $SS$-convexity condition. Recall that a \emph{qualitative game} 
    is a collection $(X_i,P_i)$, where $X_i$ is the strategy set of player $i$ and $P_i: X=\prod_{i\in I}
    X_i \rightarrow 2^{X_i}$ is the preference correspondence of player $i$. A point $x^*$ of $X$ is 
    called an \emph{equilibrium} of the game $(X_i,P_i)$ if $P_i(x^*)=\emptyset $ for each $i\in I$.
    
    \begin{theorem}
    Let ${\mathcal E}=(X_i, P_i)$ be a qualitative game such that for each $i\in I$,
    \begin{itemize}
    \item[(c1)] $X_i$ is a nonempty compact subset of a topological space;
    
    \item[(c2)] for each finite subset $\{x_{i0}, x_{i1}, \cdots, x_{in_i}\}$ of $X_i$, there exists 
    a continuous mapping $\varphi_{n_i}:\Delta_{n_i}\rightarrow X_i$ satisfying for any subset $J_i$ 
    of $\{0, 1, \cdots, n_i\}$, $x\in \bigcap_{j\in J_i}P_i^{-1}(x_{ij})$ implies 
    $\varphi_{n_i}(\Lambda_i)\in P_i(x)$ for any 
    \[
    \Lambda_i\in \{(\lambda_{i0},\lambda_{i1},\cdots, 
    \lambda_{in_i})\in \Delta_{n_i}: \lambda_{ij}=0, \mbox{ for all }j\notin J_i\};
    \] 
    and
    
    \item[(c3)] $P_i$ has open lower section.
    \end{itemize}
    Then there exists an equilibrium of ${\mathcal E}$.
    \end{theorem}
    
    \begin{proof} 
    For each $i\in I$, define a correspondence 
    \[
    A_i: X=\prod_{i\in I}X_i\rightarrow 2^{X_i}
    \]
    by setting $A_i(x)=X_i$. Then $\Gamma=(X_i,A_i,P_i)$ is an abstract economy satisfying all conditions 
    in Theorem \ref{thm:main}. Consequently, there exists an $x^*\in X$ such that $A_i(x^*)\cap 
    P_i(x^*)=\emptyset$ for each $i\in I$. Since $P_i(x^*)\subseteq X_i=A_i(x^*)$, it follows 
    $P_i(x^*)=\emptyset $ for each $i\in I$. 
    \end{proof}


    \end{document}